\documentclass[12pt]{amsart}
\usepackage{amssymb}
\usepackage{amsfonts}
\usepackage{latexsym}
\usepackage{amscd}
\usepackage[mathscr]{euscript}
\usepackage{xy} \xyoption{all}

\newtheorem{definition}{Definition}[section]
\newtheorem{theorem}[definition]{Theorem}
\newtheorem{lemma}[definition]{Lemma}
\newtheorem{corollary}[definition]{Corollary}

\newtheorem{proposition}[definition]{Proposition}

\usepackage{setspace}

\def\Z{\mathbb Z}
\def\K{\mathbb K}

\begin{document}

\title[Simplicity of Lie algebras of Leavitt algebras]{On the simplicity of Lie algebras associated to Leavitt algebras}

\author{Gene Abrams}
\address{Department of Mathematics, University of Colorado,
Colorado Springs CO 80933 U.S.A.} \email{abrams@math.uccs.edu}
 \author{Darren Funk-Neubauer}
 \address{Department of Mathematics, Colorado State University,
Pueblo CO  81001  U.S.A.} \email{darren.funkneubauer@colostate-pueblo.edu}

\thanks{The first  author is partially supported by the U.S. National Security Agency under grant number H89230-09-1-0066. } 

\begin{abstract}
For any field $\K$ and integer $n\geq 2$ we consider the Leavitt algebra $L_\K(n)$; for any integer $d\geq 1$ we form the matrix ring $S = M_d(L_\K(n))$.  $S$ is an associative algebra, but we view $S$ as a Lie algebra using the bracket $[a,b]=ab-ba$ for $a,b \in S$.  We denote this Lie algebra as $S^-$, and consider its Lie subalgebra $[S^-,S^-]$.   In our main result, we show that $[S^-,S^-]$ is a simple Lie algebra if and only if char$(\K)$ divides $n-1$ and char$(\K)$ does not divide $d$.   In particular, when $d=1$ we get that $[L_\K(n)^-,L_\K(n)^-]$ is a simple Lie algebra if and only if char$(\K)$ divides $n-1$.
\end{abstract}

\maketitle

\date{}

\maketitle

\section{Introduction}

\doublespacing

In the seminal paper \cite{Leavitt1}, Leavitt constructed and investigated properties of associative rings $R$ which lack the invariant basis property; that is, rings for which two free left $R$-modules of different finite rank are isomorphic.    In particular, for each integer $n\geq 2$ and each field $\K$, Leavitt constructed the (now-so-called) {\it Leavitt $\K$-algebra of order $n$}, denoted in this paper by $L_\K(n)$.  The algebra $R=L_\K(n)$ has the property that ${}_RR \cong {}_RR^{n}$ as left $R$-modules.

The study of Leavitt algebras has enjoyed a recent resurgence.  For instance, the Leavitt algebras were re-discovered  in \cite{McCrimmon}, where they arise as quotients of the {\it deep matrix algebras}.   In \cite{Kennedy2} the study of Lie algebras associated to deep matrix algebras and these quotients is taken up; the genesis of our article is based in part on some of the ideas presented in \cite{Kennedy2}.  (We describe the connection between Leavitt algebras and deep matrix algebras at the end of this section.)    Next,  matrix rings over Leavitt algebras were investigated in \cite{AAP}, where isomorphisms between such algebras were shown to lead to explicit isomorphisms between important classes of C$^*$-algebras.  Finally,  the Leavitt algebras have come to be viewed as specific (important) cases of the {\it Leavitt path algebras} (see e.g. \cite{AAr1}).

Throughout this article we let $\K$ denote any field and $1_{\K}$ denote the multiplicative identity of $\K$.
We begin by formally defining the algebras of interest.

\begin{definition}
\rm
\label{def:cohnandleavitt}
 Let $n\geq 2$ be an integer, and let $\K$ be any field.
 \begin{enumerate}
\item The {\it Cohn $\K$-algebra of order $n$}, denoted $C_\K(n)$, is the unital associative $\K$-algebra with generators $x_i, y_i$ $(1 \leq i \leq n)$ and relations $y_i x_j = \delta_{i,j} 1_{C_\K(n)}$ $(1 \leq i, j \leq n)$ where $1_{C_\K(n)}$ denotes the multiplicative identity of $C_\K(n)$.
\item Let $M$ denote the two sided ideal of $C_\K(n)$ generated by $1_{C_\K(n)}-\sum_{i=1}^{n}x_i y_i$.    The {\it Leavitt $\K$-algebra of order $n$}, denoted $L_\K(n)$, is the quotient algebra $C_\K(n) / M$.
\end{enumerate}
\end{definition}

Additionally, we let  $-: C_{\K}(n) \rightarrow L_{\K}(n)$ denote the canonical quotient map, so that, for $z\in C_{\K}(n)$, $z + M \in L_{\K}(n)$ will often be denoted by $\overline{z}$.

We note that $L_\K(n)$ may also be defined as the $\K$-algebra generated by the same generators and relations as those of $C_\K(n)$, plus the additional relation $1_{C_\K(n)} = \sum_{i=1}^{n} x_i y_i$.  We choose to present the above quotient definition of $L_\K(n)$ to make some of our results more transparent.  In particular, note that $1_{L_\K(n)} = \overline{1_{C_\K(n)}} = 1_{C_\K(n)}+M$.

Although the Cohn and Leavitt algebras are associative algebras, our main interest in this paper is with Lie algebras.  Any associative $\K$-algebra $A$ can be turned into a Lie algebra over $\K$ in the usual way:   for $a,b \in A$ we define the Lie bracket $[a,b]  = ab -ba$.  We denote this Lie algebra structure on $A$ by $A^-$.
 We  recall here a few basic definitions concerning Lie algebras (for more information see \cite{Humph}).
 Let $\mathfrak{L}$ denote a Lie algebra over $\K$.  $[\mathfrak{L},\mathfrak{L}]$ denotes the Lie subalgebra of $\mathfrak{L}$ consisting of all finite sums of elements of the form $[a,b]$ where $a,b \in \mathfrak{L}$.   A subset $\mathfrak{I}$ of $\mathfrak{L}$ is an {\it ideal} of $\mathfrak{L}$ if $\mathfrak{I}$ is a $\K$-subspace of the vector space $\mathfrak{L}$ and $[\mathfrak{L}, \mathfrak{I}] \subseteq \mathfrak{I}$.  The Lie algebra  $\mathfrak{L}$ is called  {\it simple} in case $[\mathfrak{L}, \mathfrak{L}] \neq \{0\}$ and the only ideals of $\mathfrak{L}$ are $\{0\}$ and $\mathfrak{L}$.

 To provide some context for our main result we present the following theorem from classical Lie theory.  For an integer $d \geq 2$ let $M_d(\K)$ denote the $d \times d$ matrix algebra over $\K$.  Then $M_d(\K)^- = \mathfrak{gl}_d(\K)$, the general linear Lie algebra over $\K$.   Let $\mathfrak{sl}_d(\K)$ denote the Lie subalgebra of $\mathfrak{gl}_d(\K)$ consisting of trace zero matrices.  It is well known that $[ \mathfrak{gl}_d(\K) , \mathfrak{gl}_d(\K)] = \mathfrak{sl}_d(\K)$ (see \cite[p. 9]{Humph}) and that $\mathfrak{sl}_d(\K)$ is simple if and only if char$(\K)$ does not divide $d$ (see \cite[Chap. 2]{Seligman}).  From these comments we obtain the following theorem. \\

\noindent
{\bf Theorem.}
For any field $\K$ and integer $d \geq 2$, the Lie algebra $[M_d(\K)^-, M_d(\K)^-]$ is simple if and only if char$(\K)$ does not divide $d$. \\

We are now in position to present our main result (Theorem \ref{thm:main}).

\medskip

{\bf Theorem.}
Let $\K$ be any field, and let $n\geq 2$, $d\geq 1$ be integers.    Then the Lie algebra $[M_d(L_\K(n))^-, M_d(L_\K(n))^-]$ is simple if and only if char$(\K)$ divides $n-1$ and char$(\K)$ does not divide $d$.

\medskip

Our main result corrects \cite[Theorem 6.10]{Kennedy2}, and extends the result both to all fields (including those of characteristic $2$), and to all-sized matrix rings.  (In \cite{Kennedy2} the Lie algebra $[L_\K(n)^-, L_\K(n)^-]$ is denoted by $\mathfrak{sla}$.)

We conclude this introductory section by establishing some notation, and presenting some results concerning the algebras $C_\K(n)$ and $L_\K(n)$.
Throughout this article $n$ denotes an integer greater than or equal to $2$, $X$  denotes the set $\{1, 2, \ldots, n \}$, and $Mon(X)$ denotes the free monoid generated by $X$ (i.e., finite words in the alphabet $X$ with multiplication  given by concatenation).  The identity element of $Mon(X)$ is the empty word $\phi$.  For nonempty $I \in Mon(X)$ we write $I = i_1 i_2 \ldots i_t$ with $t \in \Z ^+$ and $i_s \in X$ $(1 \leq s \leq t)$.   We define $I^{rev} = i_t i_{t-1} \ldots i_1$ and $\phi ^{rev} = \phi$.

$Mon(X)$ is equipped with a partial order defined by setting
$$I \leq J \ \mbox{in case there exists} \ K \in Mon(X) \ \mbox{with} \ J=IK.$$
For $I, J \in Mon(X)$ we write $I \sim J$ if either $I \leq J$ or $J \leq I$.  We write $I \nsim J$ if $I \nleq J$ and $J \nleq I$

For a nonempty $I \in Mon(X)$ with $I=i_1i_2 \ldots i_t$ we denote by  $x_{I}$ the element $x_{i_1} x_{i_2} \ldots x_{i_t}$ of $C_\K(n)$; similarly, we denote by  $y_{I}$ the element  $y_{i_1} y_{i_2} \ldots y_{i_t}$ of $C_\K(n)$.  We let $x_{\phi}$ and $y_{\phi}$ both denote the element $1_{C_\K(n)}$ of $C_\K(n)$.   It is clear  that for all $I, J \in Mon(X)$ we have $x_{IJ} = x_{I} x_{J}$ and $y_{IJ} = y_{I} y_{J}$ in $C_\K(n)$.  Also, for all $I \in Mon(X)$ we have $y_{I} \, x_{I^{rev}} = 1_{C_\K(n)}$.

\begin{lemma}\label{cohnmult}

{\rm (i)} \  The set $\{ \, x_{I} \,  y_{J} \, | \, I, J \in Mon(X) \,  \}$ is a basis for $C_\K(n)$.

{\rm (ii)} \ For all $I, J, K, L \in Mon(X)$ the product of basis elements in $C_\K(n)$ is given by \\ \\
$( x_{I} \, y_{J} ) ( x_{K} \, y_{L} ) = \left\{ \begin{array}{cl}
 0 & \mbox{if } \, J^{rev} \nsim K \\ x_{IM} \, y_{L} &  \mbox{if } \, J^{rev} \leq K \mbox{ with }  \, K=J^{rev}M \mbox{ for some } \, M \in Mon(X) \\ x_{I} \, y_{N^{rev} L} & \mbox{if } \, K \leq J^{rev} \mbox{ with }  \, J^{rev}= KN \mbox{ for some } \, N \in Mon(X). \end{array} \right.$

\end{lemma}

\smallskip

\noindent
{\it Proof:}  This is immediate from Definition \ref{def:cohnandleavitt} and the partial order on $Mon(X)$.
\hfill $\Box$ \\

The following result regarding the Leavitt algebra $L_\K(n)$ will be used later (see Lemma \ref{L(n)dimprops} and Corollary \ref{matrixdimcorollary}) to establish that the Lie algebras corresponding to the Leavitt algebras are in fact nontrivial.

\begin{lemma}\label{xsubIindependence}
Let $\{I_v \  | \  1\leq v \leq V \}$ be a set of distinct elements of $Mon(X)$, for some positive integer $V$.  Then the set $\{\overline{x_{I_v}} \  | \  1\leq v \leq V \}$ is linearly independent in $L_{\K}(n)$.

\end{lemma}

\noindent
{\it Proof:}   This follows directly from  \cite[Lemma 1.1]{SilesMolina}, a result which applies more generally to all Leavitt path algebras.   (Alternatively, the result may be directly established in $L_\K(n)$  by utilizing the natural $\Z$-grading on $L_\K(n)$ which arises by setting $deg(x_i)=1$ and $deg(y_i)=-1$ for $1\leq i \leq n$.)
\hfill $\Box$ \\

We close this section with a brief discussion of  the deep matrix algebras of finite type and their connection to the algebras in this paper (see \cite{Kennedy2} and \cite{McCrimmon} for further details).    The {\it deep matrix algebra}  $\mathcal{DM}(X, \K)$ is the unital associative $\K$-algebra with basis elements $\mathfrak{e}(h,k)$ indexed by pairs $h,k \in Mon(X)$.  The multiplication of these basis elements is defined in \cite{Kennedy2} in a fashion similar to the multiplication described in Lemma \ref{cohnmult}. For $|X|=n$, the deep matrix algebra of finite type $\mathcal{DM}(X, \K)$ is isomorphic to the Cohn algebra $C_{\K}(n)$, via the $\K$-linear extension of the map  $\psi:   \mathfrak{e}(I, J) \mapsto  x_{I} \, y_{J^{rev}}$.  Furthermore, in  \cite{Kennedy2} the author defines $\mathcal{M}$ to be the two sided ideal of $\mathcal{DM}(X, \K)$ generated by the element $\mathfrak{e}(\phi, \phi) - \sum_{x \in X} \mathfrak{e}(x,x)$.  Observe that $\psi (\mathcal{M}) = M$ (where $M$ is the ideal of $C_\K(n)$ presented in  Definition \ref{def:cohnandleavitt}(ii)), and so $\mathcal{DM}(X,\K) / \mathcal{M}$ and $L_{\K}(n)$ are isomorphic as $\K$-algebras.

\bigskip

\section{The trace functions $T$, $\tau$, and $\tau_d$}

In this section we introduce certain linear transformations $T: C_\K(n) \rightarrow \K$,  $\tau: L_\K(n) \rightarrow \K$, and $\tau_d: M_d(L_\K(n)) \rightarrow \K$  which will be used to establish one direction of our main result, Theorem \ref{thm:main}.

\begin{definition}
\label{def:T}
{\rm  Let $T: C_\K(n) \rightarrow \K$ denote the unique $\K$-linear transformation determined by the basis assignments:
\begin{eqnarray*}
T ( x_{I} \, y_{J} ) = \left\{ \begin{array}{cl}
0 & \mbox{if } \, I \neq J^{rev}, \\ 1_{\K} &  \mbox{if } \, I = J^{rev}. \end{array} \right.
\end{eqnarray*}
In particular, $T(1_{C_\K(n)}) = T(x_{\phi} \, y_{\phi}) = 1_{\K}$. }
\end{definition}

\begin{lemma}
\label{thm:cancel}
 For all $A, B, C \in Mon(X)$ and $i \in X$ we have
\begin{enumerate}
\item $T(x_{AB} \, y_{CA^{rev}}) = T(x_{B} \, y_{C})$.
\item $T(x_{A^{rev}B} \, y_{CA} ) = T(x_{B} \, y_{C})$.
\item $T(x_{Ai} \, y_{iB}) = T(x_{A} \, y_{B})$.
\end{enumerate}
\end{lemma}

\noindent
{\it Proof:}
This follows immediately from Definition \ref{def:T} and the operation {\it rev} in $Mon(X).$
\hfill $\Box$ \\

The next result yields that  $T:C_\K(n) \rightarrow \K$ satisfies a property analogous to a property of the standard trace function on matrix algebras over $\K$.

\begin{proposition}
\label{thm:traceprod}
 $T(c \, c') = T(c' \, c)$ for all $c, c' \in C_\K(n)$.
\end{proposition}

\noindent
{\it Proof:}  By Lemma \ref{cohnmult}(i) and since $T$ is linear it suffices to show that for all $I, J, K, L \in Mon(X)$
\begin{eqnarray}
\label{traceprod1}
T( \, (x_{I} \, y_{J} ) ( x_{K} \, y_{L}) \, ) = T( \, (x_{K} \, y_{L} ) ( x_{I} \, y_{J}) \, ).
\end{eqnarray}
The proof is broken into four cases depending on whether or not $J^{rev} \sim K$ and whether or not $L^{rev} \sim I$.  In the three cases:  (Case 1) $J^{rev} \nsim K$ and $L^{rev} \nsim I$;  \ (Case 2)  $J^{rev} \nsim K$ and $L^{rev} \sim I$; \ and (Case 3) $J^{rev} \sim K$ and $L^{rev} \nsim I$, the result follows directly from Lemma \ref{cohnmult}(ii) and the definition of $\sim$ in $Mon(X)$.

For Case 4,  suppose $J^{rev} \sim K$ and $L^{rev} \sim I$.    We break Case 4 into four subcases depending on whether $J^{rev} \leq K$ or $K \leq J^{rev}$ and whether $L^{rev} \leq I$ or $I \leq L^{rev}$.  In Subcase (i) we suppose $J^{rev} \leq K$  and $I \leq L^{rev}$.  Then by the definition of $\leq$ we get that $( x_{I} \, y_{J} ) ( x_{K} \, y_{L} ) = ( x_{K} \, y_{L} ) ( x_{I} \, y_{J} ) $, so (\ref{traceprod1}) holds.  Subcase (ii), in which  $K \leq J^{rev}$ and  $L^{rev} \leq I$, is established similar to Subcase (i).   In Subcase (iii), we suppose $J^{rev} \leq K$ and $L^{rev} \leq I$. Then there exists $K', I' \in Mon(X)$ such that $K = J^{rev} K'$ and $I=L^{rev}I'$.  So we have $( x_{I} \, y_{J} ) ( x_{K} \, y_{L} ) = x_{L^{rev}} \, x_{I'} \, y_{J}\, x_{J^{rev}} \, x_{K'} \, y_{L} \ = x_{L^{rev} I' K'} \, y_{\phi L}$.  From this and Lemma \ref{thm:cancel}(ii) we have
\begin{eqnarray}
\label{traceprod6}
T( \, (x_{I} \, y_{J} ) ( x_{K} \, y_{L}) \, ) = T( \, x_{I'K'} \, y_{\phi} \, ) = \left\{ \begin{array}{cl}
0 & \mbox{if } \, K' \neq \phi \,\, \mbox{or } \, I' \neq \phi, \\ 1_{\K} &  \mbox{if } \, K' = \phi \,\, \mbox{and } \, I' = \phi. \end{array} \right.
\end{eqnarray}
Similarly we have $( x_{K} \, y_{L} ) ( x_{I} \, y_{J} ) = x_{J^{rev}} \, x_{K'} \, y_{L}\, x_{L^{rev}} \, x_{I'} \, y_{J} = x_{J^{rev} K' I'} \, y_{\phi J}$.  From this and Lemma \ref{thm:cancel}(ii) we have
\begin{eqnarray}
\label{traceprod7}
T( \, (x_{K} \, y_{L} ) ( x_{I} \, y_{J}) \, ) = T( \, x_{K'I'} \, y_{\phi} \, ) = \left\{ \begin{array}{cl}
0 & \mbox{if } \, K' \neq \phi \,\, \mbox{or } \, I' \neq \phi, \\ 1_{\K} &  \mbox{if } \, K' = \phi \,\, \mbox{and } \, I' = \phi. \end{array} \right.
\end{eqnarray}
Combining (\ref{traceprod6}) and (\ref{traceprod7}) we obtain (\ref{traceprod1}). \ Subcase (iv), in which $K \leq J^{rev}$ and $I \leq L^{rev}$, is established similar to Subcase (iii).
\hfill $\Box$ \\

We now establish an important property of $T$, and subsequently use it  to define the linear transformation $\tau : L_\K(n) \rightarrow \K$.   While all of our previous results have been independent of the field $\K$, the next result imposes a relationship between char$(\K)$ and $n$.

\begin{lemma}
\label{thm:Tm=0}
Assume char$(\K)$ divides $n-1$.  Then $T(m) = 0$ for all $m \in M$.
\end{lemma}

\noindent
{\it Proof:}  By  Lemma \ref{cohnmult}(i)  we have
 $$M = {\rm span}_{\K}  \{ \, x_{I} \,  y_{J} \, (1_{C_\K(n)} - \sum_{i=1}^{n} x_i y_i ) \, x_{K} \, y_{L} \, | \, I, J, K, L \in Mon(X) \,  \}.$$
Thus it suffices to show that for  all $I, J, K, L \in Mon(X)$
\begin{eqnarray}
\label{Tm=01}
T( \, x_{I} \,  y_{J} \, (1_{C_\K(n)} - \sum_{i=1}^{n} x_i y_i ) \, x_{K} \, y_{L} \, ) = 0.
\end{eqnarray}

  Case 1:  Suppose $J \neq \phi$.  We use  Lemma \ref{cohnmult} to get,  for $1 \leq i \leq n$,
\begin{eqnarray*}
x_{I} \, y_{J} \, x_{i} \, y_{i}  = \left\{ \begin{array}{cl}
x_{I} \, y_{J} & \mbox{if $i$ is the first letter of the word $J^{rev}$} , \\ 0 &  \mbox{otherwise}. \end{array} \right.
\end{eqnarray*}
From this we conclude that
$$x_{I} \,  y_{J} \, (1_{C_\K(n)} - \sum_{i=1}^{n} x_i y_i ) \, x_{K} \, y_{L} = (x_{I} \,  y_{J} - \sum_{i=1}^{n} x_{I} \, x_{J} \, x_i y_i ) \, x_{K} \, y_{L} = 0,$$
 and so we obtain $(\ref{Tm=01})$.

  Case 2: Suppose $K \neq \phi$.  The result is established in a manner analogous to that used in Case 1.

   Case 3:  Suppose $J = K = \phi$.  Recall that $x_{\phi}=y_{\phi} = 1_{C_\K(n)}$, and observe that
   $$x_{I} \,  y_{\phi} \, (1_{C_\K(n)} - \sum_{i=1}^{n} x_i y_i ) \, x_{\phi} \, y_{L} = x_{I} \,  y_{L} \, - \, \sum_{i=1}^{n} x_{Ii} \, y_{iL}.$$
     From this and Lemma \ref{thm:cancel}(iii) we have
\begin{eqnarray*}
T(\, x_{I} \,  y_{\phi} \, (1_{C_\K(n)} - \sum_{i=1}^{n} x_i y_i ) \, x_{\phi} \, y_{L} \, ) = T( \, x_{I} \, y_{L} \, ) - \sum_{i=1}^{n} T( \, x_{I} \, y_{L} \, ) = \left\{ \begin{array}{cl}
0 & \mbox{if } I \neq L^{rev}, \\ 1_{\K} - n 1_{\K} & \mbox{if } I = L^{rev}. \end{array} \right.
\end{eqnarray*}
Since char$(\K)$ divides $n-1$ we have $1_\K - n 1_{\K} = -(n-1)1_{\K} = 0$ in $\K$.  Combining the previous two sentences we obtain (\ref{Tm=01}).
\hfill  $\Box$ \\

\begin{definition}
\rm
\label{def:tau}
Assume char$(\K)$ divides $n-1$.   Recall that $L_\K(n) = C_\K(n)/M.$    We define
$$\tau:  L_\K(n) \rightarrow \K$$ as the $\K$-linear transformation given by $\tau (c + M) = T(c)$ for $c \in C_\K(n)$.
\end{definition}

We note that an application of Lemma \ref{thm:Tm=0} is required to ensure that $\tau$ is well defined.
In particular,
$$\tau(1_{L_\K(n)}) = \tau(1_{C_\K(n)} + M) = T(1_{C_\K(n)}) = T(x_{\phi} y_{\phi}) = 1_{\K}.$$

\begin{proposition}
\label{thm:tauproduct}
Assume char$(\K)$ divides $n-1$.  Then the linear transformation \\ $\tau : L_\K(n) \rightarrow \K$   has the following properties:
\begin{enumerate}
\item $\tau(l \, l') = \tau(l' \, l)$ for all $l,l' \in L_\K(n).$
\item $\tau(1_{L_\K(n)}) = 1_\K \neq 0$ in $\K$.
\end{enumerate}
\end{proposition}

\noindent
{\it Proof:}
Both properties follow immediately from the definition of $\tau$.
\hfill $\Box$ \\

We conclude this section by defining natural extensions of $\tau$ to matrix rings over $L_\K(n)$.     We say that an associative $\K$-algebra $A$ {\it admits the trace function} $tr$ in case  there exists a $\K$-linear functional $tr: A \rightarrow \K$ for which $tr(aa') = tr(a'a)$ for all $a,a' \in A$.

\begin{lemma}
\label{traceonmatriceslemma}
 Let $A$ be an associative $\K$-algebra which admits the trace function $tr$.  Let $d$ be any positive integer.  Then the matrix algebra $M_d(A)$ admits the trace function $tr_d$ defined by setting
$$tr_d(B) = \sum_{i=1}^d tr(B_{i,i})$$
for each $B \in M_d(A)$.
\end{lemma}

\noindent
{\it Proof:}   The proof that $tr_d$ is a $\K$-linear functional is routine.   The fact that $tr_d(BC) = tr_d(CB)$ for all $B,C \in M_d(A)$ follows directly from the definition of matrix multiplication and the fact that $tr(aa') = tr(a'a)$ for all $a,a' \in A$.
\hfill $\Box$ \\

\begin{corollary}
\label{traceonmatrices}
Assume char$(\K)$ divides $n-1$.   Then the matrix algebra $M_d(L_\K(n))$ admits a trace function $\tau_d$ by setting
$$\tau_d(B) = \sum_{i=1}^d \tau(B_{i,i})$$
for each $B \in M_d(L_\K(n))$.

\end{corollary}

\noindent
{\it Proof:}
This follows immediately from Lemma \ref{traceonmatriceslemma} together with Proposition \ref{thm:tauproduct}(i).
\hfill $\Box$ \\

\section{The main result}\label{MainResultSection}

In this section we prove our main result, Theorem \ref{thm:main}.  We start by establishing some ``nontriviality" results for the Lie algebras under consideration.

\begin{lemma}
\label{L(n)dimprops}
Let $\K$ be any field, and $n\geq 2$.
\begin{enumerate}
\item $[L_\K(n)^-, L_\K(n)^-]$ is infinite dimensional over $\K$.
\item  $L_\K(n)^-$ is infinite dimensional over $\K$.
\item Let $W$ denote $[L_\K(n)^-, L_\K(n)^-]$.  Then $[W,W] \neq \{0\}$.
\end{enumerate}
\end{lemma}

\noindent
{\it Proof:}
%Let $-: C_{\K}(n) \rightarrow L_{\K}(n)$ denote the canonical quotient map.  \\
(i): Using Lemma \ref{xsubIindependence}
%\ref{def:cohnandleavitt}
 and the definition of $[-,-]$, we see that the set \\ $\{ \, [ \overline{x_1}, \overline{x_2}^j ] \, | \, j \in \Z, j\geq 1 \}$ is $\K$-linearly independent,  and the result follows.

(ii): Follows immediately from (i) and the fact  that $[L_\K(n)^-, L_\K(n)^-] \subseteq L_\K(n)^-$.

 (iii): Again using Lemma \ref{xsubIindependence} and the definition of $[-,-]$,
 %\ref{def:cohnandleavitt}
  the element $[ \, [ \overline{x_1}, \overline{x_2} ] \, , \, [ \overline{x_1}, \overline{x_2}^2] \, ]$ of $[W,W]$ is nonzero, and the result follows.
\hfill $\Box$ \\

\begin{corollary}\label{matrixdimcorollary}
Let $\K$ be any field,  $n\geq 2$, and $d\geq 1$.  Let $S$ denote the matrix algebra $M_d(L_\K(n))$.
\begin{enumerate}
\item $[S^-, S^-]$ is infinite dimensional over $\K$.
\item  $S^-$ is infinite dimensional over $\K$.
\item Let $W$ denote $[S^-, S^-]$.  Then $[W,W] \neq \{0\}$.
\end{enumerate}
\end{corollary}

\noindent
{\it Proof:}
 All three statements follow directly by analyzing the elements in the $(1,1)$ entry of $M_d(L_\K(n))$, and applying Lemma \ref{L(n)dimprops}.
\hfill $\Box$ \\

 For an associative algebra $A$ we let $Z(A)$ denote the center of $A$, that is,
  $$Z(A) = \{ \, a \in A \, | \, aa' = a'a \  \mbox{for all} \ a' \in A \, \}.$$
We will utilize the following three Theorems A, B, and C to achieve our main result.  The first of these Theorems is due to Herstein.
\medskip

\noindent
{\bf Theorem A.} (\cite[Theorem 1.13]{Herstein})
\label{thm:herstein}
Let $A$ be a simple associative $\K$-algebra.  Assume either that  char$(A) \neq 2$, or that $A$ is not $4$-dimensional over $Z(A)$.   Then $U \subseteq Z(A)$ for any proper  Lie ideal $U$ of the Lie algebra $[A^-, A^-]$.

\medskip

We will apply Theorem A to the associative ring $M_d(L_\K(n))$.  Therefore we will need the following theorem, which  was established by Leavitt in the early 1960's.

\medskip

\noindent
{\bf Theorem B.}  (\cite[Theorem 2]{Leavitt2})   For any field $\K$ and integer $n\geq 2$, $L_\K(n)$ is a simple $\K$-algebra.  Thus, using a standard ring-theoretic result, for any field $\K$, integer $n\geq 2$, and integer $d\geq 1$, the matrix ring  $M_d(L_\K(n))$ is a simple $\K$-algebra.

\medskip

(In \cite[Theorem 4.1]{FaulknerMcCrimmon} and \cite[Theorem 3]{Kennedy1} Theorem B was rediscovered in the context of deep matrix algebras.  We note that the  the simplicity of  $L_\K(n)$ does not depend on the field $\K$.)
When applying Theorem A to $M_d(L_\K(n))$ the following result will be useful.

\medskip

\noindent
{\bf Theorem C.} (\cite[Theorem 4.2]{ArandaCrow} or \cite[Theorem 15]{Kennedy1})
\label{thm:center}
The center of $L_\K(n)$ is one dimensional; that is,
$Z(L_\K(n)) = \K 1_{L_\K(n)} = \{ \, \alpha 1_{L_\K(n)} \, | \, \alpha \in \K \, \}.$
Thus, using a standard ring-theoretic result, the center of $M_d(L_\K(n))$ is one dimensional for every integer $d\geq 1$; that is, $Z(M_d(L_\K(n))) = \K 1_{M_d(L_\K(n))} = \{ \, \alpha 1_{M_d(L_\K(n))} \, | \, \alpha \in \K \, \}.$

\medskip

We combine the previous three theorems to get the following result, which will be the main ingredient in the proof of Theorem \ref{thm:main}.

\begin{proposition}
\label{thm:keyproposition}
Let $\K$ be any field, and let $n\geq 2$,  $d\geq 1$ be integers.  Then $1_{M_d(L_\K(n))} \notin [M_d(L_\K(n))^-, M_d(L_\K(n))^-]$ if and only if  $[M_d(L_\K(n))^-, M_d(L_\K(n))^-]$ is a simple Lie algebra.
\end{proposition}

\noindent
{\it Proof:} Let $W$ denote the Lie algebra $[M_d(L_\K(n))^-, M_d(L_\K(n))^-]$.

{\bf $(\Longleftarrow)$:} \,Assume $1_{M_d(L_\K(n))} \in W$.  We prove $W$ is not simple.  By hypothesis $\K 1_{M_d(L_\K(n))} \subseteq W$.  Clearly, $\K 1_{M_d(L_\K(n))} \neq \{0\}$ and by Corollary  \ref{matrixdimcorollary}(i) $\K 1_{M_d(L_\K(n))} \neq W$.  Also, $\K 1_{M_d(L_\K(n))}$ is a $\K$-subspace of $W$ and $[W, \K 1_{M_d(L_\K(n))}] = \{0\} \subseteq \K 1_{M_d(L_\K(n))}$.  Thus, $\K 1_{M_d(L_\K(n))}$ is a nonzero, proper Lie ideal of $W$ and so $W$ is not simple.

{\bf $(\Longrightarrow)$:} \,Assume $1_{M_d(L_\K(n))} \notin W$.  We prove $W$ is simple.  By Corollary \ref{matrixdimcorollary}(iii) we have $[W,W] \neq \{0\}$.  Thus it suffices to show that the only ideals of $W$ are $\{0\}$ and $W$.  Let $U$ denote a proper Lie ideal of $W$.  We show $U=\{0\}$.  By Theorem B,  $M_d(L_\K(n))$ is a simple ring.  From Corollary \ref{matrixdimcorollary}(ii) we find $M_d(L_\K(n))$ is not $4$-dimensional over $Z(M_d(L_\K(n)))$.  Thus, by Theorems A and C we get $U \subseteq \K 1_{M_d(L_\K(n))}$.  Therefore, since $1_{M_d(L_\K(n))} \notin W$ we have $U \subseteq \K 1_{M_d(L_\K(n))} \cap W = \{0\}$ and so $U=\{0\}$.
 \hfill $\Box$ \\

We now restate and prove our main result.

 \begin{theorem}\label{thm:main}
 Let $\K$ be any field, and let $n\geq 2$, $d\geq 1$ be integers.    Then the Lie algebra $[M_d(L_\K(n))^-, M_d(L_\K(n))^-]$ is simple if and only if char$(\K)$ divides $n-1$ and char$(\K)$ does not divide $d$.
 \end{theorem}

{\it Proof.}   \ \ {\bf $(\Longleftarrow)$:}   We prove the contrapositive. So we assume $[M_d(L_\K(n))^-, M_d(L_\K(n))^-]$ is not simple, and that  char$(\K)$ divides $n-1$; we seek to prove that char$(\K)$  divides $d$.    By Proposition \ref{thm:keyproposition} we have $1_{M_d(L_\K(n))} \in  [M_d(L_\K(n))^-, M_d(L_\K(n))^-]$,    so there exist a positive integer $r$, and $A_i,A'_i \in M_d(L_\K(n))^-$ $(1 \leq i \leq r)$ such that
\begin{eqnarray*}
 1_{M_d(L_\K(n))} = \sum_{i=1}^{r} [ \, A_i, A'_i \, ] =\sum_{i=1}^{r}  A_i A'_i - A'_i A_i .
\end{eqnarray*}
Since char$(\K)$ divides $n-1$, Corollary \ref{traceonmatrices} applies to yield a trace function $\tau_d: M_d(L_\K(n)) \rightarrow \K$.  By the displayed equation and the definition of trace function we get $\tau_d(1_{M_d(L_\K(n))})=0$.  On the other hand, by the definition of $\tau_d$ we have  $\tau_d(1_{M_d(L_\K(n))}) = d\cdot \tau(1_{L_\K(n)})$, while by Proposition \ref{thm:tauproduct}(ii) we have  $\tau( 1_{L_\K(n)}) \neq 0$ in $\K$.  Thus we get that $d=0$ in $\K$, so that char$(\K)$  divides $d$ as desired. \\

{\bf $(\Longrightarrow)$:} \, Assume char$(\K)$ does not divide $n-1$, or that char$(\K)$ divides $d$.  We prove in either case that the Lie algebra $[M_d(L_\K(n))^-, M_d(L_\K(n))^-]$ is not simple. By Proposition \ref{thm:keyproposition} it suffices to show $1_{M_d(L_\K(n))} \in [M_d(L_\K(n))^-, M_d(L_\K(n))^-].$

 Case 1:   Assume char$(\K)$ does not divide $n-1$.  Let $-: C_\K(n) \rightarrow L_\K(n)$ denote the canonical quotient map, that is, $\overline {c} = c+M$ for all $c \in C_\K(n)$.  For each $1\leq i,j \leq d$ and each $\ell \in L_\K(n)$ let $\ell \epsilon_{i,j}$ denote the element of $M_d(L_\K(n))$ which is $\ell$ in the $(i,j)$ entry, and zero otherwise.  Recall by Definition \ref{def:cohnandleavitt} that  we have
$\overline{1_{C_\K(n)}} = \sum_{i=1}^{n} \, \overline{x_i} \, \overline{y_i}$ in $L_\K(n)$.   Moreover, in $L_\K(n)$ we have \begin{align*}
\sum_{i=1}^{n} \, [ \, \overline{y_i}, \overline{x_i} \, ]
&=   \sum_{i=1}^{n} \, ( \, \overline{y_ix_i} - \overline {x_i} \, \overline{y_i} \, ) \\
&=   \sum_{i=1}^{n} \, \overline{1_{C_\K(n)}} - \sum_{i=1}^{n} \, \overline{x_i} \, \overline{y_i} \\
&= n \, \overline{1_{C_\K(n)}} - \overline{1_{C_\K(n)}}  \\
&= (n-1) \, 1_{L_\K(n)}.
\end{align*}
Using this, we now compute in $M_d(L_\K(n))$:
\begin{align*}
\sum_{j=1}^d \sum_{i=1}^{n} \, [ \, \overline{y_i}\epsilon_{j,j}, \overline{x_i}\epsilon_{j,j} \, ]
&=  \sum_{j=1}^d ( \sum_{i=1}^n [ \, \overline{y_i}, \overline{x_i} \, ] \, ) \epsilon_{j,j} \\
&=   \sum_{j=1}^{d} \, ((n-1) \, 1_{L_\K(n)})\epsilon_{j,j} \\
&= (n-1) \, 1_{M_d(L_\K(n))}.
\end{align*}
Using the hypothesis that  char$(\K)$ does not divide $n-1$, we multiply both sides of the above equation by $(n-1)^{-1}$ and conclude that $1_{M_d(L_\K(n))} \in [M_d(L_\K(n))^-, M_d(L_\K(n))^-]$.

\smallskip

Case 2:  Now assume char$(\K)$ divides $d$.   We note that for each $1\leq j\leq d-1$ we have
\begin{eqnarray}
\label{matrixbracket}
[1_{L_\K(n)}\epsilon_{j,j+1},1_{L_\K(n)}\epsilon_{j+1,j}] = 1_{L_\K(n)}\epsilon_{j,j} - 1_{L_\K(n)}\epsilon_{j+1,j+1} \ .
\end{eqnarray}
We multiply equation (\ref{matrixbracket}) by $j$ for each $1\leq j \leq d-1$ and then sum the results; this is easily seen to yield
$$
\sum_{j=1}^{d-1} \ j\cdot [1_{L_\K(n)}\epsilon_{j,j+1},1_{L_\K(n)}\epsilon_{j+1,j}] \  = \  \begin{pmatrix}1_{L_\K(n)}&0& &0&0 \\
        0&1_{L_\K(n)}& &0&0 \\
        & &\ddots &  &\vdots &  \\
             0&0&&1_{L_\K(n)}&0 \\
             0&0&\cdots &0&-(d-1)1_{L_\K(n)} \end{pmatrix}.
$$
Since char$(\K)$ divides $d$ we have that $-(d-1)1_{L_\K(n)} = 1_{L_\K(n)}$.  So the displayed matrix is $1_{M_d(L_\K(n))}$, so that $1_{M_d(L_\K(n))} \in [M_d(L_\K(n))^-, M_d(L_\K(n))^-]$, thus establishing the result.
\hfill $\Box$ \\

\bigskip

We record the $d=1$ case of Theorem \ref{thm:main} as a result of interest in its own right.

\begin{corollary}
Let $\K$ be any field, and let $n\geq 2$ be an integer.    Then the Lie algebra $[L_\K(n)^-, L_\K(n)^-]$ is simple if and only if char$(\K)$ divides $n-1$.
\end{corollary}

We conclude with the following observation.  Suppose $E$ is a finite directed graph, and let $L_\K(E)$ denote the Leavitt path algebra of $E$ with coefficients in $\K$ (see e.g. \cite{AAr1}).  The Leavitt path algebras include as specific examples not only the algebras $L_\K(n)$, but the matrix rings $M_d(\K)$ and  $M_d(L_\K(n))$ as well.  It would be interesting to find necessary and sufficient conditions on  $\K$ and $E$ so that the Lie algebra $[L_\K(E)^-,L_\K(E)^-]$ is simple.

\bigskip

\begin{center}
ACKNOWLEDGMENT
\end{center}
The authors are grateful to the referee for a thorough report on the original version of this article, and, in particular, for suggesting a computation which made the proof of Theorem \ref{thm:main} possible.

\singlespacing

\end{document}